\documentclass[11pt]{article}

\usepackage{hyperref}
\hypersetup{
    colorlinks=true,
    linkcolor=blue,
    filecolor=magenta,      
    urlcolor=blue,
}

  \usepackage{geometry}
\usepackage{amsmath,amssymb,amsthm}
\usepackage{graphics,epsfig,calc}
\usepackage{amsmath}

\usepackage{latexsym,epsfig,bm,amssymb}
\usepackage{xcolor}
\usepackage{amsthm,mathrsfs}

\usepackage{mathptmx}

\DeclareSymbolFont{AMSb}{U}{msb}m{n}
\DeclareSymbolFontAlphabet{\mathbb}{AMSb}

 \newcommand{\beqn}{\begin{eqnarray}}
 \newcommand{\eeqn}{\end{eqnarray}}
 \newcommand{\be}{\begin{equation}}
 \newcommand{\ee}{\end{equation}}
  
  \newcommand{\bcor}{\begin{corollary}}
 \newcommand{\ecor}{\end{corollary}}

 \newcommand{\ba}{\begin{array}}
 \newcommand{\ea}{\end{array}}

 \newcommand{\pa}{\partial}
 
  \newcommand{\ci}{\cite}
 
 \newcommand{\la}{\label}
 
 \newcommand{\rRe}{{\rm Re\5}}
 \newcommand{\fr}{\frac}

 \newcommand{\supp}{{\rm supp\5}}

\newcommand{\ov}{\overline}
\newcommand{\rot}{{\rm rot\5}}
\newcommand{\dv}{{\rm div\5}}

\newcommand{\cA}{{\cal A}}
\newcommand{\bA}{{\bf A}}
\newcommand{\bba}{{\bf a}}

\newcommand{\bn}{{\bf n}}

\newcommand{\cE}{{\cal E}}

\newcommand{\bj}{{\bf j}}

\newcommand{\cH}{{\cal H}}

\newcommand{\cX}{{\mathbb X}}

\newcommand{\mmB}{{\mathbb B}}

\newcommand{\ve}{\varepsilon}
\newcommand{\vp}{\varphi}
\newcommand{\De}{\Delta}
\newcommand{\de}{\delta}

\newcommand{\al}{\alpha}
\newcommand{\ga}{\gamma}
\newcommand{\Ga}{\Gamma}
\newcommand{\vka}{\varkappa}

\newcommand{\si}{\sigma}

\newcommand{\Om}{\Omega}
\newcommand{\na}{\nabla}

\newcommand{\lam}{\lambda}
\newcommand{\Lam}{\Lambda}

\newcommand{\bPi}{{\bm{\Pi}}}

\newcommand{\5}{{\hspace{0.5mm}}}

\newcommand{\R}{\mathbb{R}}
\newcommand{\Z}{\mathbb{Z}}

\renewcommand{\theequation}{\thesection.\arabic{equation}}
\renewcommand{\thesection}{\arabic{section}}
\renewcommand{\thesubsection}{\arabic{section}.\arabic{subsection}}
\newtheorem{theorem}{Theorem}[section]

\renewcommand{\thetheorem}{\arabic{section}.\arabic{theorem}}
\newtheorem{defin}[theorem]{Definition}

\newtheorem{lemma}[theorem]{Lemma}
\newtheorem{remark}[theorem]{Remark}

\newtheorem{corollary}[theorem]{Corollary}
\newtheorem{pro}[theorem]{Proposition}

\newcommand{\bp}{\begin{pro}}
\newcommand{\ep}{\end{pro}}
\newcommand{\bl}{\begin{lemma}}
\newcommand{\el}{\end{lemma}}
\newcommand{\bd}{\begin{defin}}
\newcommand{\ed}{\end{defin}}
\newcommand{\br}{\begin{remark}}
\newcommand{\er}{\end{remark}}

\newcommand{\Ho}{\stackrel{{\!\!\!\scriptsize o}}{H^1}}

\newcommand{\bce}{\begin{center}}
\newcommand{\ece}{\end{center}}

\begin{document}

\begin{center}

{\huge On absorbing set for 3D  Maxwell--Schr\"odinger  
\bigskip

damped driven equations in bounded region}
\bigskip\bigskip


 {\Large A. I. Komech\footnote{The research supported by the Austrian Science Fund (FWF) under Grant No. P28152-N35.}
 }
  \medskip
 \\
{\it
\centerline{Faculty of Mathematics of Vienna University}
 }
 
{\it
\centerline{Institute for Transmission Information Problems of RAS, Moscow, Russia}
 }
 
 {\it
\centerline{Mechanics-Mathematics Department, Moscow State University}
 }

 \centerline{alexander.komech@univie.ac.at}

\bigskip\bigskip

\end{center}

\begin{abstract}
We consider the 3D damped driven Maxwell--Schr\"odinger equations
in a bounded region under suitable boundary conditions.  
We establish  new a priori estimates, which provide  the existence of 
global finite energy weak solutions and 
bounded absorbing set. The proofs rely on the Sobolev type
 estimates for magnetic Schr\"odinger
operator.

  \end{abstract}
  
\tableofcontents

  \section{Introduction}

The  Maxwell--Schr\"odinger, 
 coupled equations
 form a fundamental dynamical system of Quantum Theory. 
These equations describe  crucial  phenomena of the {\it matter-radiation interaction} which is in the center of Quantum Theory and its applications.
In particular, in the applications to the design and optimal control of 
 quantum high-frequency electronic devices: laser, maser, 
 klystron, magnetron,
 traveling wave tube, synchrotron, electron microscope, and others.
 The importance of these questions was pointed out  in early paper by Kapitza \ci{K1963}.
Thus, a rigorous investigation of the 
long-time asymptotics for solutions of these equations is indispensable 
for physical applications. 

However, the 
mathematical theory of these nonlinear evolutionary 
equations is currently in an initial stage. Respectively, 
 applications of these equations rely on the perturbation
theory which cannot provide the long-time behaviour of solutions  
of these equations. 
On the other hand, almost all applications require to know the 
long-time behaviour.

Till now the design and control of quantum devices  uses mainly the quasiclassical approximation, which treats electrons as classical particles. For example,  in the fundamental monograph \ci{G2011}, the words `quantum' and `Schr\"odinger equation' were not even mentioned. The most important exception is the study of the laser and maser action, based on 1D
coupled Maxwell--Schr\"odinger  equations \ci{N1973,SSL1978}. However, these studies are not rigorous, which was one of our motivations for the present investigation.

\section{Damped driven Maxwell--Schr\"odinger  equations}
We consider
the  coupled  damped driven  Maxwell--Schr\"odinger equations (MS) 
in a bounded domain $V\subset\R^3$ with a smooth boundary $\Ga:=\pa V$.
We choose the units where $e=-1$ and $m=c=\hbar=1$. Then
in the 
Coulomb gauge $\dv\bA(x,t)\equiv 0$ the equations read
(cf.\ci{K2013, K2019phys, GNS1995, BT2009})
\be\la{MSdd}
\left\{
\ba{rcl}
\ddot \bA (x, t)\!\!&\!\!=\!\!&\!\! \De \bA (x, t) - \si \dot\bA (x, t) 
+P\bj(\cdot,t),\quad \De A^0(x,t):=-\rho(x,t)
\\\\
i\dot \psi (t)\!\!&\!\!=\!\!&\!\!(1-i\ve))H(t)\psi(t)-i\ga E(t) \psi(t)
\ea\right|,\qquad x\in V.
\ee
Here 
 $\si>0$ is the electrical conductance of the medium, 
$\ve,\ga>0$ are the absorption coefficients, and 
 $P$ denotes the orthogonal projection onto free-divergent vector fields  
from the Hilbert space $L^2(V)\otimes\R^3$.
Further,
$
E(t):=\langle\psi(t), H(t)\psi(t)\rangle,
$
where
$H(t)$ is the Schr\"odinger operator
\be \la {H}
H(t):=\fr12
D^2(t)  + \phi (x)+A^0(x,t),\qquad D(t):=-i\na+\bA(x,t)+\bA_p(x,t).
\ee
Here $\bA_p  (x, t)$ is an external `pumping potential', and
 $\phi ( x)$ stands for a  static external potential (in the case of an atom,  $\phi(x)$ is the nucleus  potential).
 Finally,
the charge and current densities are expressed
in the wave function and the Maxwell potentials as
\be \la {rji}
\rho (x, t) = |\psi(x,t)|^2, \qquad \bj (x, t) = 
\rRe [\ov{\psi(x,t)}D(t)\psi(x,t)].
\ee

\br\la{rgEi}
{\rm
We introduce 
in the Schr\"odinger equation of the system (\ref{MSdd})
the novel specific nonlinear damping term $-i\ga E(t) \psi(t)$ which plays the key role
in our approach.

}
\er

\section{Boundary conditions}
We choose the boundary conditions 
modelling  ideally conducting 
diamagnetic materials (like cooper, silver, gold, etc).
In such materials the electric and magnetic field should vanish
as well as the charge and current surface densities.
Hence, the 
  tangential component of the electric field 
  $$E(x,t)=-\dot \cA(x,t)-\na A^0(x,t),\qquad \cA(x,t):=\bA (x, t)+\bA_p (x, t)$$
vanishes on the  boundary $\Ga=\pa V$
as well as the normal component of the magnetic field $B(x,t)=\rot A(x,t)$. 
More precisely, we assume that  
\be \la {BCAi}
A^0(x,t)=0,\qquad
\bn (x)\times\dot\cA (x, t) = 0,\quad
\bn (x) \cdot \rot \5 \5 \cA (x, t) = 0, \qquad x \in \Ga,\,\,\,t>0,
\ee
where
$ \bn (x) $ is the outward normal to the resonator boundary at the point $ x \in \Ga $.
We slightly reinforce the middle boundary conditions assuming 
\be\la{BCA}
\bn(x)\times\cA (x, t)=0,\qquad x \in \Ga,\,\,\,t>0.
\ee
Then the second conditions of (\ref{BCAi}) hold by differentiation. 
Moreover, then 
the third condition
follows from (\ref{BCA}) in local orthogonal coordinates. 
Indeed,
let a point $y \in \pa V$ and $x_3=0$ on the tangent plane $T_y\Ga$.
Then (\ref{BCA}) means that 
\be\la{An}\cA_1(y,t)=\cA_2(y,t)=0\Leftrightarrow
\cA(y,t)=C(y)\bn(y),\qquad y\in\Ga.
\ee
These identities imply that 
\be\la{An2}
\pa_k\cA_j(y,t)=0,\qquad y\in\Ga, \quad k,j=1,2,
\ee
if  $\cA\in C^1(\ov V)\otimes\R^3$.
In particular,  $\pa_1\cA_2(y,t)-\pa_2\cA_1(y,t)=0$ which implies
the last  condition of (\ref{BCAi}).

Finally, for the electronic wave function we assume the Dirichlet boundary condition 
\be \la {BCp}
\psi(x,t)=0,\qquad x\in\pa V,\,\,\,t>0,
\ee
which ensures the absence of electronic current on the  boundary:
$\bj (x, t) =0$ for $x\in\pa V$.

\section{Hamiltonian structure}

The Hamiltonian functional is defined by
\be\la{encs2}
\cH(\bA,\bPi,\psi,t)=
\fr12[c^2\Vert  \bPi\Vert^2
+\Vert \rot\bA\Vert^2]+
\fr 12\langle\psi, H_0(t)\psi\rangle,
\qquad H_0(t):=\fr12D_0^2(t)  + \phi (x)+\fr12A^0(x).
\ee
Here 
$D_0(t):=-i\na+\bA(x)+\bA_p(x,t)$
and
$A^0(x):=(-\De)^{-1}\rho(\cdot)=(-\De)^{-1}|\psi(\cdot)|^2$,
where 
 $(-\De)^{-1}$ is specified with the Dirichlet boundary conditions 
for $A^0(x,t)$ from 
(\ref{BCAi}).
The system (\ref{MSdd}) 
under the boundary conditions (\ref{BCAi}),   (\ref{BCp})
can be formally written as
\be\la{MSddH}
 \dot\bA(t)=\cH_\bPi,\quad
 \dot{\bPi}(t)=-\cH_\bA- \si \bPi (x, t),
 \quad
i\dot\psi (t)=(1-i\ve)
\cH_\psi
-i\ga  E(t) \psi(t).
\ee

\section{Comments on previous results}

$\bullet$ The Maxwell--Schr\"odinger
system of type  (\ref{MSdd}) in a bounded region was not considered previously.
Such system was considered in \ci{GNS1995,BT2009}
for the case 
of the infinite space
 $V=\R^d$  
 with $d=1,2,3$ for
$\si=\ve=\ga=0$, $c=1$ and zero pumping $\bA_p(x,t)\equiv 0$:
\be \la {MS}
\left\{\ba{rcl}
\ddot \bA (x, t) \!\!&\!\!=\!\!&\!\! \De \bA (x, t) 
+P\bj (\cdot, t),\quad  \De A^0(x,t)= -\rho (x, t)
\\
\\
i\dot \psi ( t)\!\!&\!\!=\!\!&\!\!
H(t)\psi(t)
\ea\right|, \,\,\,x \in \R^d.
\ee
For this system the existence  of global solutions 
for all finite energy initial states  
was proved for the first time by Guo, Nakamitsu and  Strauss \ci {GNS1995}.
Their approach  relies on application of
 the Gagliardo--Nirenberg interpolation inequality. The uniqueness of the solution was not proved.

 The complete result on the well-posedness in the energy space 
was established 
by Bejenaru and Tataru \ci{BT2009} providing strong a priori estimates.
The methods \ci{BT2009} rely on  microlocal analysis of pseudodifferential operators with ``rough
  symbols''. In particular, these methods provide Lemma 11 of \ci{BT2009}:
 For each $0 \le s \le 2$ the operator $1-[\na-i\bA]^2$ is a diffeomorphism $H^s(\R^3)\to H^{s-2}
 (\R^3)$
which depends continuously on $\bA \in H^1(\R^3)$.
This lemma is a refinement  of   Proposition A.I of \ci{GV1981}.
 \smallskip\\
   $\bullet$ {\bf Dissipative autonomous  evolutionary PDEs.}
  The theory of attractors and long-time behaviour of solutions of
  such equations
originated in the works of Ball, Foias, Hale, Henry, Temam,  and 
was developed further 
by  Babin and Vishik, Chepyzhov, Haraux, Ilyin, Miranville,  Pata, Zelik,   and others for the Navier-Stokes,  reaction-diffusion, Ginzburg--Landau, damped wave and nonlinear  Schr\"odinger,  and sine-Gordon equations  \cite{BV1992}--\ci{KKddS2020}.
    \smallskip\\
 $\bullet$ {\bf Hamiltonian autonomous evolutionary PDEs.}
 My team has a long time experience working with the theory of global attractors
for nonlinear {Hamiltonian} evolutionary  PDEs.
I initiated this theory 
 in 1990. The 
  theory
was inspired by  Bohr's  transitions between quantum stationary states and 
resulted in
more than 50 papers
including  joint  papers  in collaboration 
with H. Spohn, V. Buslaev and others.
The global attraction to a compact attractor was established for 
a list of nonlinear Hamiltonian PDEs,
 see the surveys
\ci{K2016, KKumn2020,S2004}.
  The proofs rely on a novel application
of subtle tools of Harmonic Analysis: the Wiener Tauberian theorem,
the Titchmarsh convolution theorem, the theory of quasimeasures, and others.

Tao established 
in \ci{T2008}
the existence of the global attractor for 
radial solutions to   nonlinear 
defocusing Schr\"odinger equation without damping in $\R^n$ with $n\ge 11$.
\smallskip\\
$\bullet$ {\bf Damped  driven 
 nonlinear wave and Ginzburg-Landau equations.}  Absorbing sets and global attractors were constructed  
 i) for damped  driven 
 nonlinear wave equations 
 by Haraux \ci{Har1981,Har1991} for almost periodic external force, see also Ghidaglia and Temam 
 \ci{GT1987}, Mora and Sol\`a-Morales \ci{MSM1987},
 Babin, Chepyzhov and  Vishik \ci{BV1992}--\ci{CV2002}, and others,
 and ii) 
 for the damped  driven Ginzburg-Landau equations  
 by Ghidaglia and Heron \ci{GH1987} (see also \ci{BV1992,T1997}) in the case of
 a bounded region $V\subset\R^n$ with $n=1,2$.
 \smallskip\\
 $\bullet$ {\bf Damped driven nonlinear Schr\"odinger equations.}
The  theory of global attractors 
was developed
by Ghidaglia  \ci{G1988} on a bounded interval $V\subset \R$, 
 by Wang \ci{W1995} on the circle $V=\R/\Z$,
by Abounouh \ci{A1993} on a bounded region $V\subset\R^2$, 
and by Lauren\c{c}ot \ci{L1995} on $V=\R^N$ with $N\le 3$. 
These results on attraction were developed in \ci{GM2009,W1995}.
In  these papers  the  pumping term   does not depend on time,
and in \ci{G1988} the pumping term is time-periodic.
The main achievement in these papers is  the 
construction  of a compact global attractor in the energy space $H^1$
relying on the Ball ideas  \ci{B2004}.

A compact global attractor in the energy space for 2D weakly damped driven nonlinear Schr\"odinger equation
with general nonlinearity   
  was constructed 
in  joint paper of the PI with E. Kopylova \ci{KKddS2020}
for   general bounded region $V\subset\R^2$ and almost periodic driving.
\smallskip\\
$\bullet$  {\bf The Maxwell--Klein--Gordon and other coupled equations.} 
For various coupled equations
the results on well-posedness, existence 
of solitary waves  and their effective dynamics  
 were
obtained by P. D'Ancona, M. Esteban, S. Klainerman, M. Machedon, 
S. Selberg, E. S\'er\'e, D. Stuart, and others, \ci{GV1981}--
\ci{S2010}.
 

 
\section{Sobolev type estimates for magnetic Schr\"odinger operator}
We denote the spaces $L^p=L^p(V)$, $H^s=H^s(V)$, 
 $\Ho=\Ho(V)$, and $\Vert\cdot\Vert$ is the norm in $L^2$.
Denote $\cX=[H^1\otimes\R^3]\oplus [L^2(\R^3)\otimes\R^3]\oplus\Ho$  the Hilbert space of  states $X=(\bA,\bPi,\psi)$ satisfying the boundary conditions (\ref{BCAi})
and $\dv \bA(x)=\dv\Pi(x)=0$ for $x\in V$.
\smallskip

We will prove below   in Lemma \ref{lLam}  the  bound  for magnetic potential
\be\la{covA}
\Vert\bA\Vert_{L^2}^2\le C\Vert \na\bA\Vert^2,\qquad \bA\in\mathbb A.
\ee
This bound holds since the Laplacian $\De$ under the boundary conditions (\ref{BCA})
is nonnegative and symmetric on a dense domain $D\subset\bA$, and  hence, it  
admits the selfadjoint extension. Finally, the spectrum is discrete and zero is not an  eigenvalue. 
In Lemma \ref{lD} we prove the equivalence of norms
 for magnetic Schr\"odinger operator
\be \la{covp}
b_1(\Vert\bA(t)\Vert_{H^1}^2) \Vert\psi\Vert_{H^1}^2\le 
[\Vert D(t)\psi\Vert   +\Vert\psi\Vert]^2\le b_2(\Vert\bA(t)\Vert_{H^1}^2) \Vert\psi\Vert_{H^1}^2, \qquad\psi\in \Ho,
\ee
where
$D(t):=i\na-\bA(x,t)$,
 $b_1(r)>0$ (respectively $b_2(r)>0$) is a  decreasing (respectively an increasing) function of $r\ge 0$.
Similarly,
\be \la{covp2}
b_1(\Vert\bA(t)\Vert_{H^1}^2) \Vert\psi\Vert_{H^2}^2\le
[\Vert H(t)\psi\Vert   +\Vert\psi\Vert]^2\le b_2(\Vert\bA(t)\Vert_{H^1}^2) \Vert\psi\Vert_{H^2}^2, \qquad\psi\in H^2\cap\Ho.
\ee
The  bounds (\ref{covp}) and (\ref{covp2}) 
extend   Lemma 11 of \ci{BT2009} to the 
case of bounded region.
For the proof of (\ref{covp}) we will  show that the difference of $D^2(t)$ with $-\De$
is a relatively compact operator. 
\smallskip

We will assume that
 the potential $\phi(x)$ is bounded, 
\be\la{V}
\sup_{x\in V} |
\phi(x)|<\infty.
\ee
 Hence, we can assume that
it is positive,
\be\la{pp}
\phi(x)\ge \vka >0,\qquad x\in V
\ee
since the potential is defined up to an additive constant.
Hence,
\be\la{ppS}
E(t)=
\langle\psi,H(t)\psi\rangle\ge \Vert D(t)\psi\Vert^2   +\vka \Vert\psi\Vert^2
+\langle\rho, (-\De)^{-1}\rho\rangle.
 \ee
 Now (\ref{covp}) implies 
 \be\la{ppS2}
 E(t)\ge \vka_1
 (\Vert\bA(t)\Vert_{H^1}^2) \Vert\psi\Vert_{H^1}^2,
 \ee
 where $\vka_1>0$ is a decreasing function.
Hence, the standard  Sobolev estimates together with (\ref{ppS}) imply from
(\ref{covA}) and (\ref{covp})
 that
\be\la{cov2}
\left\{\ba{rcll}
\Vert\bA\Vert_{L^p}^2&\le& C\sum_k\Vert \na_k\bA\Vert^2,& \bA\in\mathbb A
\\\\
\Vert\psi\Vert_{L^p}^2&\le& b(\Vert\bA(t)\Vert_{H^1}) E(t)
,&\psi\in \Ho
\ea\right| ,\qquad p\in[2,6].
 \ee
 Here the  last bound 
extends   Lemma 11 of \ci{BT2009} to the 
case of bounded region.

 \section{A priori estimates}
 First, 
 we obtain a priori estimates for sufficiently
 smooth solutions
$(\bA(x,t),\bPi(x,t),\psi(x,t))$ of the Maxwell--Schr\"odinger system (\ref{MSdd}), where all functions are $C^\infty(\R^4)$.
We plan to get rid of this smoothness assumption and establish the same 
estimates for all finite energy solutions.
We will assume for the pumping field $\bA_p(x,t)$ that it is 
{\it almost periodic} and 
\be\la{pum}
\sup_{x\in V,t\in\R}
[|\bA_p(x,t)|+|\na \bA_p(x,t)|+|\dot \bA_p(x,t)|]<\infty.
\ee
Differentiating the charge $Q(t):=\Vert\psi(t)\Vert^2$ we have
\beqn\la{enc2Q}
\dot Q(t)\!\!\!&\!\!\!=\!\!\!&\!\!\!\langle\dot\psi(t),\psi(t)\rangle+\langle\psi(t),\dot\psi(t)\rangle
 \nonumber\\
 \nonumber\\
\!\!\!&\!\!\!=\!\!\!&\!\!\! 
\langle (-i-\ve)H(t)\psi(t),\psi(t)\rangle+\langle \psi(t),(-i-\ve)H(t)\psi(t)\rangle
-2\ga
 E(t)
 Q(t)
 \nonumber\\
 \nonumber\\
\!\!\!&\!\!\!=\!\!\!&\! -2\ve
 E(t)
-2\ga
 E(t)
 Q(t)
 \le
 -2\ve\vka Q(t)
 -2\ga\vka  Q^2(t),
\eeqn
where we used (\ref{ppS}).
Hence,
\be\la{enc3s}
Q(t)\le Q(0).
\ee
Now 
differentiating the energy $\cE(t):=\cH(\bPi(t),\bA(t),\psi(t),t)$ 
and using (\ref{MSdd}) and (\ref{ppS}),  (\ref{pum}), we get
\beqn
\la{enc2s}
\!\!\!\!\!\!\dot\cE(t)\!\!\!&\!\!=\!\!\!&\!\!\langle\cH_\bA,\dot\bA\rangle+\langle\cH_\bPi,\dot\bPi\rangle+\langle\cH_\psi,\dot\psi\rangle+\cH_t
\nonumber\\
\nonumber\\
\!\!\!\!\!\!\!\!\!&\!\!=\!\!\!&\!\!
\langle\cH_\bA,\cH_\bPi\rangle\!+\!\langle\cH_\bPi,-\cH_\bA-\si\bPi\!+\!
\langle\cH_\psi,(-i-\ve)\cH_\psi\!-\!\ga
 E(t) \psi(t)-\langle D(t)\psi,\dot \bA_p\psi\rangle
\nonumber\\\nonumber\\
\!\!\!\!\!\!\!\!\!&\!\!\le\!\!\!&\!\!
-\si\Vert\bPi(t)\Vert^2
\!-\!\ve\langle H(t)\psi(t), H(t)\psi(t)\rangle
\!-\!\ga E^2(t) 
+C_p E(t)+C_p\Vert D(t)\psi(t)\Vert\Vert\psi(t)\Vert
\!\le\! C_1<\infty,
\eeqn
since $\cH_\psi=H(t)\psi(t)$.
Hence,
 (\ref{H}) and
  (\ref{covp}), and (\ref{covp2}) imply a priori estimate
\be\la{enc4s}
\Vert\na\bA(t)\Vert^2\!+\!\Vert\bPi(t)\Vert^2\!+\!
\Vert \psi(t)\Vert_{H^1}^2\!+\!\ve \!\int_0^t\!
b_1(\Vert\bA(s)\Vert_{H^1})
\Vert \psi(s)\Vert_{H^2}^2ds
\!\le\! C(t\!+\!1),\,\, \Vert\psi(t)\Vert\le C\!<\!\infty,\quad  t>0.
\ee

\section{Absorbing set}
The estimates  (\ref{enc4s})
are insufficient to prove the existence of a bounded  absorbing set.
 We will 
follow the ideas of   Haraux \ci{Har1981} 
(see also \ci{BV1992})
introducing the functional 
\be\la{HM}
\Phi(\bA,\bPi,\psi,t)=\cH(\bA,\bPi,\psi,t)+\eta\langle\bPi,\bA\rangle
\ee
with  a small $\eta>0$.
 Differentiating $\Phi(t):=\Phi(\bA(t),\bPi,\psi(t),t)$ 
 and using (\ref{enc2s}), we obtain
 \beqn\la{HMP}
\dot\Phi(t)\!\!\!&\!\!\!=\!\!\!&\!\!\!\dot\cE(t)+\eta\langle\bPi(t),\bPi(t)\rangle
+\eta\langle\dot\bPi(t),\bA(t)\rangle
\nonumber\\
\nonumber\\
\!\!\!&\!\!\!\le \!\!\!&\!\!\!
-\si\Vert\bPi(t)\Vert^2
\!-\!\ve\langle H(t)\psi(t), H(t)\psi(t)\rangle
\!-\!\ga E^2(t)+C_p E(t)+\eta\langle\bPi(t),\bPi(t)\rangle
\nonumber\\
\nonumber\\
\!\!\!&\!\!\!&\!\!\!
+\eta\langle
\De \bA (t) - \si\bPi ( t) 
+\bj(t),\bA(t)
\rangle
\eeqn
The most problematic term 
\be\la{Aj}
\langle\bj(t),\bA(t)\rangle=-\rRe\langle\ov\psi(t) D(t)\psi(t),\bA(t)\rangle
\ee
can be estimated using the  Sobolev-type estimates
 (\ref{cov2}):
\beqn\la{Aj2}
|\langle\bj(t),\bA(t)\rangle|&\le&
 C\Vert\bA(t)\Vert_{L^6}\Vert\psi(t)\Vert_{L^3}\cdot\Vert D(t)\psi(t)\Vert
 \le C_1 
 \Vert\na\bA(t)\Vert
  E^{1/2}(t)
 \Vert D(t)\psi(t)\Vert
\nonumber\\
\nonumber\\
&\le &
\de\Vert\na\bA(t)\Vert^2+\fr {C_2}\de  E^2(t),
\eeqn
where the last inequality holds by (\ref{ppS}).
For the remaining terms similar estimates  follow from the first estimate (\ref{covA}) 
and from (\ref{pum}):
\be\la{rt}
|\langle\bPi(t),\bA(t)\rangle|\le \de 
\Vert\na\bA(t)\Vert^2+\fr 1\de \Vert\bPi(t)\Vert^2.
\ee
Now (\ref{HMP}) implies
that for any $\de>0$
\be\la{HMP2}
\dot\Phi(t)\le
-\eta(1-3\de)\Vert\na\bA(t)\Vert^2
-(\si-\eta-\fr\eta\de)\Vert\bPi(t)\Vert^2
-(\ga-C_2\fr\eta\de)  E^2(t)
-\ve\Vert H(t)\psi(t)\Vert^2
+C_p.
\ee
It remains to  choose $\de,\eta>0$ such that 
\be\la{al} 
\min(\eta(1-3\de), \si-\eta-\fr\eta\de, \ga-C_2\fr\eta\de)>0.
\ee
Then  (\ref{HMP2}) and (\ref{encs2})   imply that
\be\la{HMP3}
\dot\Phi(t)\le
-\al \cH(t)-\ve\Vert H(t)\psi(t)\Vert^2+C_p
,\qquad t>0,
\ee
where 
$\al>0$ and  $C_p\in\R$ do not depend  on the solution. 
However,  
(\ref{covA}) implies that
for sufficiently small $\eta>0$
we have 
\be\la{PcH}
c\cH(\bA,\bPi,\psi,t)\le \Phi(\bA,\bPi,\psi,t)\le C\cH(\bA,\bPi,\psi,t)
\ee
with $c,C>0$. Hence, ( \ref{HMP3}) and ( \ref{covp2})  imply
that for small $\eta>0$
\be\la{HMP4}
\dot\Phi(t)
\le-\beta \Phi(t)-\ve\Vert H(t)\psi(t)\Vert^2+ C_p
,\qquad t>0,
\ee
where 
$\beta>0$ and  $ C_p\in\R$ do not depend  on the solution. 
Now the integration yields
\be\la{HMP5}
\Phi(t)+
\ve\int_0^t e^{-\beta(t-s)}\Vert H(s)\psi(s)\Vert^2ds\le    \Phi(0)e^{-\beta t}+\fr { C_p}\beta,\qquad t>0.
\ee
Hence,
 (\ref{encs2}) and  
 (\ref{ppS}),
 (\ref{PcH}) imply that
 for sufficiently small $\eta>0$ a priori estimate (\ref{enc4s}) refines to
\be\la{b1}
\Vert\na\bA(t)\Vert^2
\!+\!
\Vert\bPi(t)\Vert^2\!+\!
\Vert D(t)\psi\Vert^2\!  +\! \Vert\psi\Vert^2
\!+\ve\!\int_0^t\! e^{-\beta(t-s)}\Vert H(s)\psi(s)\Vert^2ds
\!\le\!  C[ \Phi(0)e^{-\beta t}\!+\!\fr { C_p}\beta],\quad  t\!>\!0.
\ee
Now  (\ref{covp}) and (\ref{covp2}) imply that
\be\la{b2}
\Vert\na\bA(t)\Vert^2\!+\!
\Vert\bPi(t)\Vert^2\!+\!
b_1(M)\Vert \psi\Vert_{H^1}^2\! +\ve  b_1(M)
\int_0^t\! e^{-\beta(t-s)}
\Vert \psi(s)\Vert_{H^2}^2
ds
\!\le\!  C[ \Phi(0)e^{-\beta t}\!+\!\fr { C_p}\beta],\quad  t\!>\!0.
\ee
where 
$$
M:=\sup_{t\ge 0} \Vert\bA(t)\Vert_{H^1}\le C  \sup_{t\ge 0} \Vert\na\bA(t)\Vert_{H^1}<\infty 
$$ by (\ref{b1}).
Let us write the system (\ref{MSdd}) as
\be\la{MSdd2X}
\dot X(t)=F(X(t),t),\qquad X(t)=(\bA(t),\bPi(t),\psi(t)).
\ee

\bcor 
The bounds  (\ref{b2}) imply for solutions $X(t)$ of  (\ref{MSdd2X}) 
\be\la{b2c}
\Vert X(t)\Vert_\cX^2
+\ve_1 
\int_0^t\! e^{-\beta(t-s)}
\Vert \psi(s)\Vert_{H^2}^2
ds
\le  C(\Phi(0)e^{-\beta t}+\fr { C_p}\beta),\qquad  t>0,
\ee
where $\ve_1,\beta>0$.
Hence, for any $R>0$ the set 
$\mmB:=\{Y\in \cX:\Vert Y\Vert_\cX^2\le C(1+\fr {C_p}\beta)\}$ absorbs the ball 
$\{Y\in \cX:\Vert Y\Vert_\cX\le R\}$ for large times $t>t_R$.

\ecor


\appendix

\protect\renewcommand{\thesection}{\Alph{section}}
\protect\renewcommand{\theequation}{\thesection.\arabic{equation}}
\protect\renewcommand{\thesubsection}{\thesection.\arabic{subsection}}
\protect\renewcommand{\thetheorem}{\Alph{section}.\arabic{theorem}}

\section{Proof of the Sobolev type estimates for magnetic Schr\"odinger operator}
Let a point $y \in \Ga:=\pa \Om$ and $x_3=0$ on the tangent plane $T_y\Ga$.
Then (\ref{An2}) together with $\dv\bA(y,t)\equiv 0 $ for $y\in\Om$ implies that
\be\la{An3}
\pa_3\bA_3(y,t)=0
\ee
if  $\bA\in C^1(\ov\Om)\otimes\R^3$.
Hence, the boundary conditions (\ref{An}), (\ref{An3})
can be written as
\be\la{BCA2}
\bA_\Vert(x)=0,\qquad \na_\bn\bA_\bn(x)=0,\qquad x\in\Ga,
\ee
where $\bA_\Vert(x)$ is the tangential to the boundary projection,
while $\bA_\bn(x)$  is the normal to the boundary projection of $\bA(x)$.
These boundary conditions 
and the Stokes formula
imply that for $\bA\in C^2(\ov\Om)\otimes\R^3$
\be\la{naA}
\langle\bA(t),\De \bA (t)\rangle=\int_\Ga \bA(x,t)\cdot\na_\bn\bA(x,t)dx-\sum_k\Vert\na_k\bA(t)\Vert^2=-\sum_k\Vert\na_k\bA(t)\Vert^2,
\ee
since the field $\bA(x)$ is orthogonal to the boundary $\Ga$ while $\na_\bn\bA(x)$ is 
parallel to the boundary at any point $x\in\Ga$ due to the boundary conditions (\ref{BCA2}).
Hence, 
the variational derivative
\be\la{An4}
D_\bA\sum_k\Vert \na_k\bA\Vert^2=-2\De\bA
\ee
if  $\bA\in C^2(\ov\Om)\otimes\R^3$ and satisfies the boundary conditions (\ref{BCA2}).

\bl\la{lLam}
The Laplacian  $\Lam=-\De$  is symmetric and nonnegative
on the dense 
domain $D_0=\cA(\Om)\cap C^\infty(\ov\Om)\otimes \R^3$
in the Hilbert space of vector fields
$X:=L^2(\Om)\otimes \R^3$,
and 
 the identity holds
\be\la{AAl}
\langle\bA,\Lam\bA\rangle=\sum_k\Vert\na_k\bA\Vert^2,\qquad \bA\in D_0.
\ee
The operator $\Lam$
admits a  selfadjoint extension with a domain $D\subset\cA$, and
\be\la{Ld}
\langle\bA,\Lam\bA\ge\de \Vert\bA\Vert^2,\qquad\bA\in D,
\ee
where $\de >0$.
\el
\begin{proof}
The Green formula implies that for any vector fields $\bA_1,\bA_2\in D_0$
\be\la{Gr}
\langle \Lam \bA_1,\bA_2\rangle-\langle \bA_1,\Lam\bA_2\rangle=-\int_\Ga[\na_\bn\bA_1(x)\cdot\bA_2(x)-
\bA_1(x)\cdot\na_\bn\bA_2(x)]dx=0
\ee
which follows 
similarly to (\ref{naA}).
Moreover,  the operator $\Lam$ is nonnegative on the domain  $D_0$
by (\ref{naA}), and 
hence, it admits a selfadjoint extension  by the Friedrichs theorem.
The identity  (\ref{AAl}) follows from  (\ref{naA}) .

Finally, $\Lam$ is an elliptic operator in the bounded region $\Om$, and 
the boundary conditions (\ref{BCA2}) satisfy the Shapiro--Lopatinski condition.
Hence,  
the spectrum of $\Lam$ is discrete.
Namely,
$\Lam-z$ is invertible for $z<0$
 and the resolvent $R(z)=(\Lam-z)^{-1}$ is a compact selfadjoint operator 
in $X$ by the elliptic theory  and the Sobolev embedding theorem.
Finally, the spectra of $\Lam-z$ and of $\Lam$ differ by the shift. 

Now to prove  (\ref{Ld})
it suffices to check that $\lam=0$ is not an eigenvalue. Indeed,
 $\Lam\bA=0$ implies $\bA\in D$ by the elliptic theory,
and hence (\ref{AAl}) implies 
\be\la{AAl0}
\langle\bA,\Lam\bA\rangle=\sum_k\Vert\na_k\bA\Vert^2=0.
\ee
Hence, $\bA(x)\equiv \bba \in\R^3$ for $x\in\Om$. Finally, the boundary conditions
(\ref{BCA2}) imply that $\bba=0$.
\end{proof}

To prove (\ref{covp}) let us rewrite it equivalently as
\be\la{covp23}
\langle \psi,\Lam\psi \rangle \le C_2(\Vert\bA\Vert_{H^1}) 
\langle \psi,L\psi \rangle,
 \qquad\psi\in \Ho,
\ee
where $\Lam :=-\De+1$ and $L:=(i\na+\bA(x))^2+1$. Now  (\ref{covp23}) follows from the next lemma with $\de<1$.

\bl\la{lD}
Let $A\in\cA$. Then
the difference $T:=L-\Lam$ admits the estimate
\be\la{T}
|\langle\psi, T\psi\rangle|\le \de \langle\psi, \Lam\psi\rangle+C_\de(\Vert\bA\Vert_{H^1}) \langle\psi, \psi\rangle, \qquad\psi\in \Ho
\ee
for any $\de>0$, where $C_\de(\cdot)$ is a continuous increasing function on $[0,\infty)$.

\el
\begin{proof}
The difference $T$ reads
\be\la{Tr}
T=2i\bA(x)\na+\bA^2(x),
\ee
where we have used that $\dv\bA(x)\equiv 0$.
Further, it suffices to prove (\ref{T}) for  $\psi\in D_0:=C_0^\infty(\Om)$
which is dense in the space $\Ho$ by its definition.
Hence, the operators $\Lam$ and $T$ can be considered as the operators on entire space
$\R^3$. 
In particular, let us extend $\bA(x)$ by zero outside $\Om$ and
 define the powers $\Lam^s$ for $s\in\R$ by the Fourier transform,
Then  (\ref{T}) can be reduced
by the substitution $\psi=\Lam^{-1/2}\vp$
to the  estimate
\be\la{T2}
|\langle\vp, \Lam^{-1/2}T\Lam^{-1/2}\vp\rangle|\le \de \Vert\vp\Vert^2+C_\de(\Vert\bA\Vert_{H^1}) 
\Vert \Lam^{-1/2}\vp\Vert^2, \qquad\vp\in L^2(\R^3).
\ee
This reduction is not equivalent, but  (\ref{T2}) implies  (\ref{T})
since every function $\psi\in D_0$ admits the representation $\psi=\Lam^{-1/2}\vp$ with 
$\vp\in L^2(\R^3)$.

In the case $\bA(x)\in C_0^\infty(\R^3)$ the operator $\Lam^{-1/2}T\Lam^{-1/2}$
is the PDO of  order $-1$, and estimates of type (\ref{T2}) in this case follows
from the interpolation inequality for the Sobolev norms.
However, the constant $C_\de$ is known to depend on  some derivatives of the symbol
of the composition. So we should 
find another arguments to
prove that this constant depends only on the norm $\Vert\bA\Vert_{H^1}(\Om)$.

For this purpose we note that  the operator $\Lam^{-1/2}$ is 
the multiplication by $\langle\xi\rangle^{-1/2}$ in the Fourier transform.
Hence, it is 
the convolution with a 
distribution $S(x)\in L^1_{\rm loc}(\R^3)$ which is a 
radial smooth function for $x\ne 0$
and asymptotically  homogeneous  at the origine
\be\la{Sh}
S(x)\sim C|x|^{-2},\qquad|x|\to 0. 
\ee
Hence, 
\be\la{Sh2}
S\in L^q_{\rm loc}(\R^3),\qquad q\in[1,\fr32).
\ee
Thus,  $\Lam^{-1/2}=S*$,  and hence,
\be\la{Tr2}
\Lam^{-1/2}T\Lam^{-1/2}=2iS*\bA(x)\na S*+S*\bA^2(x)S*,
\ee
Let us estimate the first term on the right hand side. The second can be bounded
similarly. 

We denote by $\na S*$  the composition  of $\na$ with $S*$. 
This composition is the bounded operator in $L^2(\R^3)$ since it
is the multiplication by $-i\xi\langle\xi\rangle^{-1/2}$
in the Fourier transform. Thus,
\be\la{Tr3}
\Vert\na S*\vp\Vert\le C\Vert\vp\Vert.
\ee
Furter, the multiplication  by $\bA\in L^6$ maps continuously $L^2(\R^3)$ into $L^{3/2}(\R^3)$
by the H\"older inequality.
\be\la{Tr4}
\Vert\bA(x)\na S*\vp\Vert_{L^{3/2}(\R^3)}\le C\Vert\vp\Vert.
\ee
Moreover,  $\supp\bA\subset\ov\Om$, and hence,
 the convolution  $S*[\bA(x)\na S*\vp]\in L^{3-\al'}(\Om)$ with sufficiently small $\al'>0$
 by
 (\ref{Sh2})
 and
  the Young theorem on the convolution.
 Similarly, 
 \be\la{TR5}
 \Vert \Lam^\al S*[\bA(x)\na S*\vp]\Vert_{L^3(\Om)}\le C\Vert\vp\Vert
 \ee
  for small $\al>0$.
 This follows by the same Young theorem since $\Lam^\al S*$ 
 for $\al<1/2$
 is 
 the operator of
 convolution with the distribution $S^\al\in L^1_{\rm loc}(\R^3)$
 which admits the asymptotics
 \be\la{Sh2al}
S^\al(x)\sim C|x|^{-2-2\al},\qquad|x|\to 0.
\ee
 This asymptotics holds since the Fourier transform of $S^\al$ equals to $\langle\xi\rangle^{-1/2+\al}$.
 
 Finally, (\ref{TR5}) means that the operator $\Lam^\al S*[\bA(x)\na S*\vp]$ is bounded
 in $L^2(\Om)$ 
  since the region $\Om$ is bounded. 
 In other words, the operator $ S*[\bA(x)\na S*\vp]$ is continuous from 
 $L^2(\Om)$ to the Sobolev space $H^{2\al}(\Om)$ with sufficiently small $\al>0$. 
 Hence, the bound of type  (\ref{T2}) for the first term on the right hand side of (\ref{Tr2}) 
 follows from the interpolation inequality for the Sobolev norms.
 The second term can be bounded similarly.
  \end{proof}



\end{document}